\documentclass[12pt,reqno]{amsart}
\usepackage{amssymb}
\usepackage{amsxtra}
\usepackage{mathrsfs}
\usepackage[all]{xy}
\newtheorem{theorem}{Theorem}[section]
\newtheorem{lemma}[theorem]{Lemma}
\newtheorem{prop}[theorem]{Proposition}
\newtheorem{corollary}[theorem]{Corollary}

\theoremstyle{definition}
\newtheorem{definition}{Definition}[section]
\theoremstyle{remark}
\newtheorem{example}{Example}[section]

\DeclareMathOperator{\Tor}{Tor}
\DeclareMathOperator{\Ext}{Ext}
\newcommand*{\ptens}[1]{\mathop{\widehat\otimes}_{#1}}

\newcommand*{\Ptens}{\mathop{\widehat\otimes}}

\newcommand*{\Tens}{\mathop{\otimes}}
\newcommand*{\Smash}{\mathop{\#}}
\newcommand*{\Psmash}{\mathop{\widehat{\#}}}

\newcommand*{\lmod}{\mbox{-}\!\mathop{\mathbf{mod}}}
\newcommand*{\rmod}{\mathop{\mathbf{mod}}\!\mbox{-}}
\newcommand*{\bimod}{\mbox{-}\!\mathop{\mathbf{mod}}\!\mbox{-}}
\newcommand*{\id}{1}

\newcommand*{\wh}{\widehat}
\newcommand*{\wt}{\widetilde}

\newcommand*{\h}{\mathbf h}

\newcommand*{\CC}{\mathbb C}

\renewcommand*{\O}{\mathscr O}

\newcommand*{\fg}{\mathfrak g}
\newcommand*{\fh}{\mathfrak h}
\newcommand*{\fr}{\mathfrak r}

\newcommand*{\cT}{\mathcal T}

\DeclareMathOperator{\Ker}{Ker}

\renewcommand*{\Im}{\mathop{\mathrm{Im}}}

\newcommand*{\eps}{\varepsilon}
\newcommand*{\lar}{\leftarrow}

\newcommand*{\xla}{\xleftarrow}
\newcommand*{\xra}{\xrightarrow}
\begin{document}
\title[Arens-Michael enveloping algebras]{Arens-Michael enveloping algebras\\
and analytic smash products}
\author{A. Yu. Pirkovskii}
\thanks{Partially supported by the RFBR grants 03-01-06392 and 02-01-00928,
and by the President of Russia grant MK-2049.2004.1.}
\date{}
\begin{abstract}
Let $\fg$ be a finite-dimensional complex Lie algebra, and
let $U(\fg)$ be its universal enveloping algebra.
We prove that if $\wh{U}(\fg)$, the Arens-Michael envelope of $U(\fg)$,
is stably flat over $U(\fg)$ (i.e., if the canonical homomorphism
$U(\fg)\to\wh{U}(\fg)$
is a localization in the sense of Taylor \cite{T2}), then $\fg$ is solvable.
To this end, given a cocommutative Hopf algebra $H$ and an $H$-module
algebra $A$, we explicitly describe the Arens-Michael envelope of
the smash product $A\# H$ as an ``analytic smash product'' of
their completions w.r.t. certain families of seminorms.
\end{abstract}
\maketitle

The Arens-Michael envelope of a complex associative algebra $A$
is defined as the completion of $A$ w.r.t. the family of all submultiplicative
seminorms on $A$. This notion (under a different name)
was introduced by Taylor \cite{T1}, and the terminology
``Arens-Michael envelope'' is due to
Helemskii \cite{X2}. An important example is the polynomial algebra
$\CC[t_1,\ldots ,t_n]$ whose Arens-Michael envelope
is isomorphic to the algebra $\O(\CC^n)$ of entire functions endowed with
the compact-open topology. Thus the Arens-Michael envelope of a
noncommutative finitely generated algebra can be viewed as
an ``algebra of noncommutative entire functions'' (cf. \cite{T2,T3}).

Given an algebra $A$, it is natural to ask to what extent homological
properties of its Arens-Michael envelope $\wh{A}$ (considered as
a topological algebra) are related to those of $A$. To handle this
problem, it is convenient to use the notion of {\em localization}.
Roughly speaking, a topological algebra
homomorphism $A\to B$ is a localization if it identifies
the category of topological $B$-modules with a full subcategory
of the category of topological $A$-modules, and if the homological
relations between $B$-modules do not change when the modules
are considered as $A$-modules. Localizations were introduced by Taylor \cite{T2}
in connection with the functional calculus problem
for several commuting Banach space operators.
A purely algebraic counterpart of this notion
was studied by Neeman and Ranicki \cite{NR}.
(Note that their terminology differs from Taylor's one;
namely, a homomorphism $A\to B$ is a localization in Taylor's sense
precisely when $B$ is stably flat over $A$ in the sense of Neeman and Ranicki.)

Thus a natural question is whether or not $\wh{A}$ is stably flat
over $A$. Taylor \cite{T2} proved that this is the case for
$A=\CC[t_1,\ldots ,t_n]$ and for $A=F_n$, the free algebra on
$n$ generators. In the case where $A=U(\fg)$, the universal enveloping algebra
of a complex Lie algebra $\fg$, Dosiev \cite{Dos_hprop} proved
that $\wh{U}(\fg)$ is stably flat over $U(\fg)$ provided $\fg$ is metabelian.
In \cite{Pir_stb} we extended this result to the case where $\fg$ admits
a positive grading. On the other hand, Taylor \cite{T2} showed that
$\wh{U}(\fg)$ is not stably flat over $U(\fg)$ if $\fg$ is semisimple.
Here we generalize this result and show that $\wh{U}(\fg)$ can be stably
flat over $U(\fg)$ only when $\fg$ is solvable.

Our approach is based on a notion of ``analytic smash product'',
which is a continuous version of the corresponding algebraic
notion \cite{Sweedler}. To prove the above-mentioned result, we first
show that for each cocommutative Hopf algebra $H$ and each $H$-module
algebra $A$ the Arens-Michael envelope of the (algebraic) smash product
$A\# H$ is isomorphic to the analytic smash product of $\wh{H}$
and the completion of $A$ w.r.t. a certain family of seminorms
determined by the action of $H$.

\section{Preliminaries}
We shall work over the complex numbers $\CC$. All associative algebras
are assumed to be unital.

By a topological algebra we mean a topological vector space $A$
together with the structure of associative algebra such that
the multiplication map $A\times A\to A$ is separately continuous.
A complete, Hausdorff, locally convex topological algebra with
jointly continuous multiplication is called a
\textit{$\Ptens$-algebra} (see \cite{T1,X1}). If $A$ is a $\Ptens$-algebra, then
the multiplication $A\times A\to A$ extends
to a linear continuous map from the completed projective tensor
product $A\Ptens A$ to $A$. In other words, a
$\Ptens$-algebra is just an
algebra in the tensor category $(\mathbf{LCS},\Ptens)$
of complete Hausdorff locally convex spaces.
The latter observation can be used to define $\Ptens$-coalgebras,
$\Ptens$-bialgebras, and Hopf $\Ptens$-algebras; see, e.g., \cite{BFGP}.

Recall that a seminorm $\|\cdot\|$ on an algebra $A$ is called
{\em submultiplicative} if $\| ab\|\le\| a\|\| b\|$ for all $a,b\in A$.
This means precisely that the corresponding unit ball $U=\{ a\in A : \| a\|\le 1\}$
is {\em idempotent}, i.e., satisfies $U^2\subset U$.
A topological algebra $A$ is said to be {\em locally $m$-convex}
if its topology can be defined by a family of submultiplicative seminorms.
Note that the multiplication in a locally $m$-convex algebra is
jointly continuous.
An {\em Arens-Michael algebra} is a complete, Hausdorff,
locally $m$-convex algebra.

The following useful lemma is due to Mitiagin, Rolewicz, and
\.Zelazko \cite{MRL}.

\begin{lemma}
Let $A$ be a locally convex algebra with topology generated
by a family $\{\|\cdot\|_\nu : \nu\in\Lambda\}$ of seminorms.
Suppose that for each $\nu\in\Lambda$ there exist $\mu\in\Lambda$ and $C>0$
such that $\| a_1\ldots a_n\|_\nu\le C^n\| a_1\|_\mu\ldots \| a_n\|_\mu$
for each $a_1,\ldots ,a_n\in A$. Then $A$ is locally $m$-convex.
\end{lemma}

\begin{corollary}
\label{cor:MRZ1}
Let $A$ be a locally convex algebra with topology generated
by a family $\{\|\cdot\|_\nu : \nu\in\Lambda\}$ of seminorms.
Suppose that for each $\nu\in\Lambda$ there exist $\mu\in\Lambda$ and $C>0$
such that $\| ab\|_\nu\le C\| a\|_\mu\| b\|_\nu$
for each $a,b\in A$. Then $A$ is locally $m$-convex.
\end{corollary}

We shall use the latter corollary in the following geometric form.

\begin{corollary}
\label{cor:MRZ2}
Let $A$ be a locally convex algebra. Suppose that $A$ has a base
$\mathscr U$ of absolutely convex $0$-neighborhoods with the property that
for each $V\in\mathscr U$ there exist $U\in\mathscr U$ and $C>0$
such that $UV\subset CV$. Then $A$ is locally $m$-convex.
\end{corollary}

Let $A$ be a topological algebra. A pair $(\wh{A},\iota_A)$ consisting of
an Arens-Michael algebra $\wh{A}$ and a continuous homomorphism
$\iota_A\colon A\to\wh{A}$ is called {\em the Arens-Michael envelope} of $A$
\cite{T1,X2}
if for each Arens-Michael algebra $B$ and for each continuous
homomorphism $\varphi\colon A\to B$ there exists a unique
continuous homomorphism $\wh{\varphi}\colon\wh{A}\to B$ making
the following diagram commutative:
\begin{equation*}
\xymatrix{
\wh{A} \ar@{-->}[r]^{\wh{\varphi}} & B \\
A \ar[u]^{\iota_A} \ar[ur]_\varphi
}
\end{equation*}
In the above situation, we say that {\em $\wh{\varphi}$ extends $\varphi$}
(though $\iota_A$ is not injective in general; see \cite{X2} or \cite{Pir_stb}
for details).

Recall (see \cite{T1} and \cite[Chap. V]{X2}) that the Arens-Michael envelope
of a topological algebra $A$ always exists and
can be obtained as the completion\footnote[1]{Here we follow the convention
that the completion of a non-Hausdorff locally convex space $E$
is defined to be the completion of the corresponding Hausdorff
space $E/\overline{\{ 0\}}$.} of $A$
w.r.t. the family of all continuous submultiplicative seminorms on $A$.
This implies, in particular, that $\iota_A\colon A\to\wh{A}$ has dense range.
Clearly, the Arens-Michael envelope is unique in the obvious sense.

Each associative algebra $A$ becomes a topological algebra
w.r.t. the finest locally convex topology. The Arens-Michael
envelope, $\wh{A}$, of the resulting topological algebra will be referred
to as the Arens-Michael envelope of $A$. That is, $\wh{A}$ is
the completion of $A$ w.r.t. the family of {\em all} submultiplicative
seminorms. Thus a neighborhood base at $0$ for the topology on $A$
inherited from $\wh{A}$ consists of all absorbing, idempotent,
absolutely convex subsets.

Here is a basic example: the Arens-Michael envelope of the polynomial
algebra $\CC[t_1,\ldots ,t_n]$ is topologically isomorphic to the
algebra $\O(\CC^n)$ of entire functions endowed with the compact-open
topology \cite{T1}. For other examples, see \cite{T1,T2,X2,Pir_stb}.

Let $\fg$ be a finite-dimensional complex Lie algebra. We may define
{\em the Arens-Michael enveloping algebra} of $\fg$ as a pair
$(\wh{U}(\fg),\iota_\fg)$ consisting of
an Arens-Michael algebra $\wh{U}(\fg)$ and a Lie algebra homomorphism
$\iota_\fg\colon\fg\to\wh{U}(\fg)$ such that for each
Arens-Michael algebra $B$ and for each Lie algebra
homomorphism $\varphi\colon\fg\to B$ there exists a unique
$\Ptens$-algebra homomorphism $\wh{\varphi}\colon\wh{U}(\fg)\to B$
such that $\wh{\varphi}\iota_\fg=\varphi$.
Clearly, $\wh{U}(\fg)$ is nothing but the Arens-Michael envelope
of $U(\fg)$, the universal enveloping algebra of $\fg$.

If $H$ is a bialgebra (resp., a Hopf algebra), then it is easy
to show that the Arens-Michael envelope $\wh{H}$ is a $\Ptens$-bialgebra
(resp., a Hopf $\Ptens$-algebra) in a natural way (for details, see \cite{Pir_stb}).

Let $A$ be a $\Ptens$-algebra.
A left \textit{$A$-$\Ptens$-module} is a complete Hausdorff
locally convex space $X$ together
with the structure of left unital $A$-module such that the map
$A\times X\to X,\; (a,x)\mapsto a\cdot x$ is jointly continuous.
As above, this means precisely that $X$ is a left $A$-module in $(\mathbf{LCS},\Ptens)$.
If $X$ and $Y$ are left $A$-$\Ptens$-modules, then
the vector space of all (continuous) $A$-module morphisms from
$X$ to $Y$ is denoted by ${_A}\h (X,Y)$. Left $A$-$\Ptens$-modules
and their (continuous) morphisms form a category denoted by $A\lmod$.
Right $A$-$\Ptens$-modules, $A$-$\Ptens$-bimodules, and their morphisms
are defined similarly. The corresponding categories are denoted
by $\rmod A$ and $A\bimod A$, respectively.

If $X$ is a right $A$-$\Ptens$-module and $Y$
is a left $A$-$\Ptens$-module, then their $A$-module tensor product
$X\ptens{A}Y$ is defined to be
the completion of the quotient $(X\Ptens Y)/N$, where $N\subset X\Ptens Y$
is the closed linear span of all elements of the form
$x\cdot a\otimes y-x\otimes a\cdot y$
($x\in X$, $y\in Y$, $a\in A$)\footnote[2]{To avoid confusion, we note that
this definition of $X\ptens{A} Y$ (due to Helemskii \cite{X2})
is different from that given by
Kiehl and Verdier \cite{KV} and Taylor \cite{T1}.
More precisely, $X\ptens{A} Y$ is the completion
of the Kiehl-Verdier-Taylor tensor product.}.
As in pure algebra, the $A$-module tensor product can be characterized
by a certain universal property (see \cite{X1} for details).

A morphism $\sigma\colon X\to Y$ of left $A$-$\Ptens$-modules
is said to be an {\em admissible epimorphism} if there exists
a linear continuous map $\tau\colon Y\to X$ such that
$\sigma\tau=\mathbf 1_Y$, i.e., if $\sigma$ is a retraction
when considered in the category of topological vector spaces.
A chain complex $X_\bullet=(X_n,d_n)$ of left $A$-$\Ptens$-modules
is called {\em admissible} if it splits as a complex
of topological vector spaces. Equivalently, $X_\bullet$ is admissible if
each $d_n$ is an admissible epimorphism of $X_{n+1}$ onto
$\Ker d_{n-1}\subset X_n$.

An $A$-module $P\in A\lmod$ is called \textit{projective} if for each
admissible epimorphism $X\to Y$ in $A\lmod$ the induced map
${_A}\h(P,X)\to {_A}\h(P,Y)$ is surjective.
Given a left $A$-$\Ptens$-module $X$, a \textit{projective resolution} of $X$
is a chain complex $P_\bullet=(P_n,d_n)_{n\ge 0}$
consisting of projective left $A$-$\Ptens$-modules $P_n$
together with a morphism $\epsilon\colon P_0\to X$
such that the augmented sequence
$$
0\xla{} X\xla{\epsilon} P_0\xla{d_0}\cdots \xla{} P_n\xla{d_n} P_{n+1}\xla{}\cdots
$$
is an admissible complex. The category $A\lmod$ has enough
projectives, i.e., each $A$-$\Ptens$-module has a projective
resolution \cite{X1}. Therefore one can define
the derived functors $\Ext$ and $\Tor$ following the general patterns
of relative homological algebra. For details, see \cite{X1}.

Similar definitions apply to right $A$-$\Ptens$-modules
and to $A$-$\Ptens$-bimodules. A projective resolution of $A$
considered as a $\Ptens$-bimodule over itself is called a {\em projective
bimodule resolution} of $A$.

Let $A$ and $B$ be $\Ptens$-algebras
and $\theta\colon A\to B$ a continuous homomorphism.
Following Taylor \cite{T2},
we say that $\theta$ is a {\em localization}
if the following conditions are satisfied:
\begin{itemize}
\item[(i)]
There exists a projective bimodule resolution
$P_\bullet\to A\to 0$
of $A$ such that the complex
\begin{equation*}
B\ptens{A} P_\bullet \ptens{A} B \to
B\ptens{A} A \ptens{A} B
\cong B\ptens{A} B \to 0
\end{equation*}
is admissible;
\item[(ii)]
The map
$B\ptens{A} B\to B,\; b_1\otimes b_2\mapsto b_1 b_2$
is a topological isomorphism.
\end{itemize}
In this situation, we say (following Neeman and Ranicki \cite{NR})
that $B$ is {\em stably flat} over $A$.

The following observation is due to Taylor \cite{T2}.

\begin{prop}
\label{prop:Tor}
Let $A\to B$ be a localization. Then for each $M\in\rmod B$
and each $N\in B\lmod$ there are natural isomorphisms
\begin{equation*}
\Tor_n^A(M,N)\cong\Tor_n^B(M,N)\qquad (n\ge 0).
\end{equation*}
\end{prop}

For later reference, let us recall a standard notation from
the theory of topological vector spaces.
Let $E$ and $F$ be locally convex spaces.
For each $0$-neighborhood $U\subset E$ and each
$0$-neighborhood $V\subset F$ let $\Gamma(U\Tens V)$
denote the absolutely convex hull of the set
\begin{equation*}
U\Tens V=\{ u\otimes v : u\in U,\; v\in V\}\subset E\Tens F.
\end{equation*}
Then all sets of the form $\Gamma(U\Tens V)$ form a base
of $0$-neighborhoods for the projective tensor product topology
on $E\Tens F$.

\section{Algebraic and analytic smash products}
Let $H$ be a bialgebra. Recall that an {\em $H$-module algebra} is
an algebra $A$ endowed with the structure of left $H$-module such that
the product $A\Tens A\to A$ and the unit map $\CC\to A$ are $H$-module
morphisms. For example, if $\fg$ is a Lie algebra acting on $A$ by derivations,
then the action $\fg\times A\to A$ extends to a map $U(\fg)\times A\to A$
making $A$ into a $U(\fg)$-module algebra. Similarly, if $G$ is a group
acting on $A$ by automorphisms, then $A$ becomes a $\CC G$-module algebra,
where $\CC G$ denotes the group algebra of $G$.

Given an $H$-module algebra $A$, the {\em smash product algebra} $A\Smash H$
is defined as follows (see, e.g., \cite{Sweedler}).
As a vector space, $A\Smash H$ is equal to $A\Tens H$.
To define multiplication,
denote by $\mu_{H,A}\colon H\Tens A\to A$ the action of $H$ on $A$,
and define $\tau\colon H\Tens A\to A\Tens H$ as the composition
\begin{equation}
\label{tau}
H\Tens A \xra{\Delta_H\otimes\id_A} H\Tens H\Tens A \xra{\id_H\otimes c_{H,A}}
H\Tens A\Tens H \xra{\mu_{H,A}\otimes\id_H} A\Tens H
\end{equation}
(Here $c_{H,A}$ denotes the flip $H\Tens A\to A\Tens H$).
Then the map
\begin{equation}
\label{smashprod}
(A\Tens H)\Tens (A\Tens H) \xra{\id_A\otimes\tau\otimes\id_H}
A\Tens A\Tens H\Tens H \xra{\mu_A\otimes\mu_H} A\Tens H
\end{equation}
is an associative multiplication on $A\Tens H$. The resulting algebra
is denoted by $A\Smash H$ and is called the {\em smash product} of $A$ with $H$.
For later reference,
note that the maps $i_1\colon A\to A\Smash H,\; i_1(a)=a\otimes 1$
and $i_2\colon H\to A\Smash H,\; i_2(h)=1\otimes h$ are algebra homomorphisms.

Similar definitions apply in the $\Ptens$-algebra case.
Namely, if $H$ is a $\Ptens$-bialgebra, then an {\em $H$-$\Ptens$-module algebra}
is a $\Ptens$-algebra $A$ together with the structure of
left $H$-$\Ptens$-module such that the product $A\Ptens A\to A$
and the unit map $\CC\to A$ are $H$-module morphisms.
We define the {\em analytic smash product} $A\Psmash H$ to be $A\Ptens H$ as a locally
convex space. By replacing $\Tens$ by $\Ptens$
in \eqref{tau} and \eqref{smashprod}, we obtain a multiplication
on $A\Psmash H$ making it into a $\Ptens$-algebra.

\begin{example}
Let $A$ be a Banach algebra, and let $G$ be a discrete group
acting on $A$ by isometric automorphisms. Then $A\Psmash\ell^1(G)$
is isomorphic to the covariance algebra $\ell^1_A(G)$ introduced
by Doplicher, Kastler, and Robinson \cite{DKR} in a more general
setting of locally compact groups.

For numerous related constructions and references, see \cite{Schwtzr}.
\end{example}

Let $E$ be a vector space, and let $\cT$ be a set of linear
operators on $E$.
\begin{definition}
We say that a seminorm $\|\cdot\|$ on $E$ is {\em $\cT$-stable} if
for each $T\in\cT$ there exists $C>0$ such that $\| Tx\|\le C\| x\|$
for each $x\in E$.
A subset $U\subset E$ is said to be {\em $\cT$-stable} if
for each $T\in\cT$ there exists $C>0$ such that $T(U)\subset CU$.
\end{definition}
Clearly, a seminorm $\|\cdot\|$ is $\cT$-stable if and only if
the unit ball $\{ x\in E : \| x\|\le 1\}$ is $\cT$-stable.

If $E$ is a left module over an associative algebra $B$, then
we say that a seminorm $\|\cdot\|$ on $E$ (resp. a subset $U\subset E$)
is {\em $B$-stable} if it is stable w.r.t. the set of operators
$\{ x\mapsto b\cdot x : b\in B\}$. Similar definitions apply in the
case where $E$ is a left module over a Lie algebra $\fg$ or a left module
over a group $G$.
Note that if a subset $M\subset B$ generates $B$ as an algebra,
then a seminorm $\|\cdot\|$ on $E$ is $B$-stable if and only if it is
$M$-stable. In particular, a seminorm on a $\fg$-module
(resp., on a $G$-module)
is $U(\fg)$-stable (resp., $\CC G$-stable) if and only if it is $\fg$-stable
(resp., $G$-stable).

\begin{definition}
Let $H$ be a bialgebra and $A$ an $H$-module algebra.
We define the {\em $H$-completion} $\wt{A}$ to be the
completion of $A$ w.r.t. the family of all $H$-stable, submultiplicative
seminorms.
\end{definition}

It is immediate from the definition that $\wt{A}$ is an Arens-Michael algebra.

\begin{prop}
Let $H$ be a bialgebra and $A$ an $H$-module algebra.
Then the action of $H$ on $A$ uniquely extends to an action of $\wh{H}$
on $\wt{A}$, so that $\wt{A}$ becomes an $\wh{H}$-$\Ptens$-module algebra.
Moreover, the smash product $\wt{A}\Psmash\wh{H}$ is an Arens-Michael algebra.
\end{prop}
\begin{proof}
Let us endow $H$ and $A$ with the topologies
inherited from $\wh{H}$ and $\wt{A}$, respectively.
In order to prove the first assertion,
it suffices to show that the action $H\times A\to A$ is jointly continuous.
Let $\|\cdot\|$ be an $H$-stable, submultiplicative seminorm on $A$.
For each $h\in H$, set $\| h\|'=\sup\{ \| h\cdot a\| : \| a\|\le 1\}$.
Since $\|\cdot\|$ is $H$-stable, it follows that $\|\cdot\|'$ is a well-defined
seminorm on $H$. We obviously have $\| h_1 h_2\|'\le \| h_1\|' \| h_2\|'$
for all $h_1,h_2\in H$, and $\| h\cdot a\| \le \| h\|' \| a\|$ for each
$h\in H,\; a\in A$.
Hence the action $H\times A\to A$ is jointly continuous, so it
uniquely extends to a jointly continuous bilinear map $\wh{H}\times \wt{A}\to \wt{A}$.
Since the canonical image of $H$ (resp., $A$) is dense in $\wh{H}$
(resp., in $\wt{A}$), it follows that $\wt{A}$ becomes an $\wh{H}$-$\Ptens$-module
algebra.

To prove that $\wt{A}\Psmash\wh{H}$ is an Arens-Michael algebra, it suffices
to show that the algebraic smash product $A\Smash H$ is locally $m$-convex w.r.t. the
projective tensor product topology. Recall that a typical $0$-neighborhood
in $A\Smash H$ has the form $\Gamma(U\Tens V)$, where
$U\subset A$ and $V\subset H$ are absorbing, absolutely convex, idempotent subsets,
and $U$ is $H$-stable. Given such $U$ and $V$, define
\begin{equation*}
W=\{ w\in H\, : \, w\cdot U\subset U\}.
\end{equation*}
It is easy to see that $W$ is absorbing, absolutely convex, and idempotent,
so it is a $0$-neighborhood in $H$. Set $V'=\Delta^{-1}(\Gamma(W\Tens V))$,
where $\Delta\colon H\to H\Tens H$ is the comultiplication on $H$.
We claim that
\begin{equation}
\label{UV'}
\Gamma(U\Tens V')\, \Gamma(U\Tens V) \subset \Gamma(U\Tens V).
\end{equation}
Indeed, take $u_1,u_2\in U,\; v_1\in V'$, and $v_2\in V$.
Then $\Delta(v_1)\in\Gamma(W\Tens V)$. Since $W\cdot U\subset U$,
i.e., $\mu_{H,A}\bigl(\Gamma(W\Tens U)\bigr)\subset U$, it follows
that
\begin{multline*}
\tau(v_1\otimes u_2) \in
(\mu_{H,A}\otimes\id_H)(\id_H\otimes c_{H,A})\bigl(\Gamma(W\Tens V\Tens U)\bigr)\\
=(\mu_{H,A}\otimes\id_H)\bigl(\Gamma(W\Tens U\Tens V)\bigr)
\subset\Gamma(U\Tens V).
\end{multline*}
Therefore,
\begin{align*}
(u_1\otimes v_1)(u_2\otimes v_2) &=
(\mu_A\otimes\mu_H)(u_1\otimes\tau(v_1\otimes u_2)\otimes v_2)\\
& \in (\mu_A\otimes\mu_H)(U\Tens \Gamma(U\Tens V)\Tens V)\subset
\Gamma(U\Tens V).
\end{align*}
This proves \eqref{UV'}. Together with Corollary \ref{cor:MRZ2}, this
implies that $A\Smash H$ is locally $m$-convex, so that $\wt{A}\Psmash\wh{H}$ is an
Arens-Michael algebra.
\end{proof}

\begin{theorem}
Let $H$ be a cocommutative Hopf algebra and $A$ an $H$-module algebra.
Then the canonical map $A\Smash H\to \wt{A}\Psmash\wh{H}$
extends to a $\Ptens$-algebra isomorphism
\begin{equation*}
(A\Smash H)\sphat\, \cong \wt{A}\Psmash\wh{H}.
\end{equation*}
\end{theorem}
\begin{proof}
Let $\varphi\colon A\Smash H\to B$ be a homomorphism to an
Arens-Michael algebra $B$. As before, we endow $A$ and $H$ with
the topologies inherited from $\wt{A}$ and $\wh{H}$, respectively.
Since the canonical image of $A\Smash H$
is dense in $\wt{A}\Psmash\wh{H}$, it suffices to show that $\varphi$
is continuous w.r.t. the projective tensor product topology on $A\Smash H$.

Define homomorphisms $\varphi_1\colon A\to B$
and $\varphi_2\colon H\to B$ by $\varphi_1(a)=\varphi(a\otimes 1)$
and $\varphi_2(h)=\varphi(1\otimes h)$. Then
\begin{equation*}
\varphi(a\otimes h)=\varphi\bigl((a\otimes 1)(1\otimes h)\bigr)=
\varphi_1(a)\varphi_2(h)
\end{equation*}
for each $a\in A,\; h\in H$. Therefore we need only prove that $\varphi_1$
and $\varphi_2$ are continuous.

Let $\|\cdot\|$ be a continuous submultiplicative seminorm on $B$.
Then the seminorms $a\mapsto \| a\|'=\| \varphi_1(a)\|\; (a\in A)$
and $h\mapsto \| h\|''=\| \varphi_2(h)\|\; (h\in H)$
are submultiplicative. This implies, in particular, that $\varphi_2$ is continuous.
To prove the continuity of $\varphi_1$, we have to show that $\|\cdot\|'$
is $H$-stable.

Let $h\in H$ be a primitive element. Then for each $a\in A$ we have
\begin{multline*}
(1\otimes h)(a\otimes 1)
=\tau(h\otimes a)
=(\mu_{H,A}\otimes\id_H)(\id_H\otimes c_{H,A})
\bigl((h\otimes 1+1\otimes h)\otimes a\bigr)\\
=(\mu_{H,A}\otimes\id_H)(h\otimes a\otimes 1+1\otimes a\otimes h)
=h\cdot a\otimes 1+a\otimes h.
\end{multline*}
Therefore,
\begin{align}
\label{prim}
\| h\cdot a\|' &= \|\varphi(h\cdot a\otimes 1)\|\notag\\
&= \|\varphi\bigl((1\otimes h)(a\otimes 1)-a\otimes h\bigr)\|\notag\\
&\le \|\varphi(1\otimes h)\|\|\varphi(a\otimes 1)\| + \|\varphi(a\otimes h)\|\notag\\
&= \|\varphi_2(h)\|\|\varphi_1(a)\|+\|\varphi_1(a)\varphi_2(h)\| \le
2C\|\varphi_1(a)\| = 2C \| a\|',
\end{align}
where $C=\|\varphi_2(h)\|$.

Now let $g\in H$ be a group-like element.
Then for each $a\in A$ we have
\begin{align*}
(1\otimes g)(a\otimes 1)
=\tau(g\otimes a)
=(\mu_{H,A}\otimes\id_H)(\id_H\otimes c_{H,A})
(g\otimes g\otimes a)\\
=(\mu_{H,A}\otimes\id_H)(g\otimes a\otimes g)
=g\cdot a\otimes g.
\end{align*}
Therefore,
\begin{align}
\label{grp}
\| g\cdot a\|' &= \|\varphi(g\cdot a\otimes 1)\|
= \|\varphi\bigl((1\otimes g)(a\otimes 1)(1\otimes g^{-1})\bigr)\|\notag\\
&\le \|\varphi(1\otimes g)\|\|\varphi(a\otimes 1)\|\|\varphi(1\otimes g^{-1})\|\notag\\
&= \|\varphi_2(g)\| \|\varphi_1(a)\| \|\varphi_2(g^{-1})\| = C \| a\|',
\end{align}
where $C=\|\varphi_2(g)\| \|\varphi_2(g^{-1})\|$.

Since $H$ is cocommutative, it is generated by primitive
and group-like elements \cite[13.1]{Sweedler}. Therefore
it follows from \eqref{prim} and \eqref{grp} that $\|\cdot\|'$ is
$H$-stable. Hence $\varphi_1$ is continuous.
In view of the above remarks,
$\varphi$ is also continuous, and so it uniquely extends to a $\Ptens$-algebra
homomorphism $\wt{A}\Psmash\wh{H}\to B$. This completes the proof.
\end{proof}

\begin{corollary}
\label{cor:Usmash}
Let $\fg$ be a Lie algebra acting on an algebra $A$ by derivations.
Then $(A\Smash U(\fg))\sphat\,\cong \wt{A}\Psmash\wh{U}(\fg)$
as $\Ptens$-algebras.
\end{corollary}

\begin{corollary}
Let $G$ be a group acting on an algebra $A$ by automorphisms.
Then $(A\Smash \CC G)\sphat\,\cong \wt{A}\Psmash\wh{\CC G}$
as $\Ptens$-algebras.
\end{corollary}

\section{The main result}

Let $H$ be a $\Ptens$-bialgebra with counit $\eps\colon H\to\CC$,
and let $A$ be an $H$-$\Ptens$-module algebra.
\begin{lemma}
\label{lemma:idepstau}
Define $\tau\colon H\Tens A\to A\Tens H$ by \eqref{tau}. Then
$(\id_A\otimes\eps)\tau=\mu_{H,A}$.
\end{lemma}
\begin{proof}
This is a direct computation:
\begin{multline*}
(\id_A\otimes\eps)\tau
=(\id_A\otimes\eps)(\mu_{H,A}\otimes\id_H)(\id_H\otimes c_{H,A})
(\Delta\otimes\id_A)\\
=\mu_{H,A} (\id_H\otimes\id_A\otimes\eps)(\id_H\otimes c_{H,A})
(\Delta\otimes\id_A)\\
=\mu_{H,A} (\id_H\otimes\eps\otimes\id_A)(\Delta\otimes\id_A)=\mu_{H,A}.
\end{multline*}
\end{proof}

\begin{lemma}
There is a unique left $A\Psmash H$-$\Ptens$-module structure on $A$ such that
\begin{equation}
\label{Amod}
(a\otimes 1)\cdot b=ab,\quad (1\otimes h)\cdot b=h\cdot b
\end{equation}
for each $a,b\in A,\; h\in H$.
\end{lemma}
\begin{proof}
Consider the map $\id_A\otimes\eps\colon A\Psmash H\to A$.
Let us prove that $\Ker(\id_A\otimes\eps)$ is a left ideal of $A\Psmash H$.
In view of the direct sum decomposition
$A\Ptens H=(A\Ptens\Ker\eps)\oplus (A\Ptens\CC 1)$,
is suffices to show that $(a_1\otimes h_1)(a_2\otimes h_2)\in\Ker(\id_A\otimes\eps)$
whenever $h_2\in\Ker\eps$. We have
\begin{multline*}
(a_1\otimes h_1)(a_2\otimes h_2)
=(\mu_A\otimes\mu_H)\bigl(a_1\otimes\tau(h_1\otimes a_2)\otimes h_2\bigr)\\
\in (\mu_A\otimes\mu_H)(A\Ptens A\Ptens H\Ptens\CC h_2)
\subset\Ker(\id_A\otimes\eps).
\end{multline*}
Therefore $\Ker(\id_A\otimes\eps)$ is a left ideal of $A\Psmash H$, so that
we can make $A$ into a left $A\Psmash H$-$\Ptens$-module
via the identification
$A=(A\Psmash H)/\Ker(\id_A\otimes\eps)$.

Now take $a,b\in A$ and $h\in H$. We have
\begin{equation*}
(a\otimes 1)\cdot b=(\id_A\otimes\eps)\bigl((a\otimes 1)(b\otimes 1)\bigr)=
(\id_A\otimes\eps)(ab\otimes 1)=ab.
\end{equation*}
On the other hand, Lemma \ref{lemma:idepstau} implies that
\begin{equation*}
(1\otimes h)\cdot b
=(\id_A\otimes\eps)\bigl((1\otimes h)(b\otimes 1)\bigr)
=(\id_A\otimes\eps)\bigl(\tau(h\otimes b)\bigr)
=h\cdot b.
\end{equation*}
Hence conditions \eqref{Amod} are satisfied.

The uniqueness readily follows from the identity $a\otimes h=(a\otimes 1)(1\otimes h)$.
\end{proof}

From now on, we endow $A$ with the left $A\Psmash H$-$\Ptens$-module
structure defined in the previous lemma.

\begin{lemma}
\label{lemma:Aproj}
Suppose that $\CC$ is a projective left $H$-$\Ptens$-module.
Then $A$ is a projective left $A\Psmash H$-module.
\end{lemma}
\begin{proof}
Since $\CC$ is projective, there exists an $H$-module morphism
$\lambda\colon\CC\to H$ such that $\eps\lambda=\id_{\CC}$.
Then the element $x_0=\lambda(1)$ satisfies $\eps(x_0)=1$ and
$hx_0=\eps(h)x_0$ for each $h\in H$. Consider the map
\begin{equation*}
\rho\colon A\to A\Psmash H,\quad \rho(a)=a\otimes x_0.
\end{equation*}
We claim that $\rho$ is a left $A\Psmash H$-module morphism.
To prove the claim, it is convenient to consider $A\Psmash H$
as a left $A$-$\Ptens$-module and as a left $H$-$\Ptens$-module
via the embeddings $i_1\colon A\to A\Psmash H$ and $i_2\colon H\to A\Psmash H$
given by $a\mapsto a\otimes 1$ and $h\mapsto 1\otimes h$,
respectively. Thus we have to show that $\rho$ is an $A$-module morphism
and an $H$-module morphism.

For each $a,b\in A$ we have
\begin{equation*}
\rho(ab)=ab\otimes x_0=(a\otimes 1)(b\otimes x_0)=a\cdot\rho(b),
\end{equation*}
so that $\rho$ is an $A$-module morphism. Further, the relation $hx_0=\eps(h)x_0$
implies that $(\id_A\otimes\mu_H)(u\otimes x_0)=(\id_A\otimes\eps)(u)\otimes x_0$
for each $u\in A\Ptens H$. Together with Lemma~\ref{lemma:idepstau}, this gives
\begin{multline*}
h\cdot\rho(a)=(1\otimes h)(a\otimes x_0)
=(\id_A\otimes\mu_H)\bigl(\tau(h\otimes a)\otimes x_0\bigr)\\
=(\id_A\otimes\eps)\bigl(\tau(h\otimes a)\bigr)\otimes x_0
=h\cdot a\otimes x_0
=\rho(h\cdot a)
\end{multline*}
for each $a\in A,\; h\in H$. Therefore $\rho$ is an $H$-module morphism and
hence an $A\Psmash H$-module morphism. Finally, since $\eps(x_0)=1$,
we see that $(\id_A\otimes\eps)\rho=\id_A$. Thus $A$ is a retract of
$A\Psmash H$ in $A\Psmash H\lmod$, so it is projective.
\end{proof}

Now let $\fg$ be a finite-dimensional complex Lie algebra.
Denote by $\fr$ the radical of $\fg$, and consider the Levi decomposition
$\fg=\fr\oplus\fh$. The action of $\fh$ on $\fr$ by commutators extends
to an action of $\fh$ on $U(\fr)$ by derivations, and there exists a canonical
isomorphism $U(\fg)\cong U(\fr)\Smash U(\fh)$ (see, e.g., \cite[1.7.11]{MR}).
Using Corollary~\ref{cor:Usmash}, we see that
\begin{equation*}
\wh{U}(\fg) \cong \bigl(U(\fr)\Smash U(\fh)\bigr)\sphat\cong
\wt{U}(\fr)\Psmash\wh{U}(\fh).
\end{equation*}

\begin{lemma}
\label{lemma:Uproj}
$\wt{U}(\fr)$ is a projective $\wh{U}(\fg)$-$\Ptens$-module.
As a corollary,
\begin{equation*}
\Tor_k^{\wh{U}(\fg)}\bigl(\CC,\wt{U}(\fr)\bigr)=0\;\text{ for each }k>0.
\end{equation*}
\end{lemma}
\begin{proof}
Since $\fh$ is semisimple, the Arens-Michael envelope
$\wh{U}(\fh)$ is isomorphic to a direct product of full matrix algebras
\cite[Corollary 7.6]{T2}. Hence each $\wh{U}(\fh)$-$\Ptens$-module is projective
\cite{T1} (see also \cite[5.28]{X1}). Now it remains to apply Lemma~\ref{lemma:Aproj}.
\end{proof}

\begin{lemma}
\label{lemma:ne0}
Suppose that $k=\dim\fh>0$. Then
$\Tor_k^{U(\fg)}\bigl(\CC,\wt{U}(\fr)\bigr)\ne 0$.
\end{lemma}
\begin{proof}
Set $A=\wt{U}(\fr)$, and recall that the groups
$H_p(\fg,A)=\Tor_p^{U(\fg)}(\CC,A)$ can be computed as
the homology groups of the standard complex $C_\cdot(\fg,A)$:
\begin{equation*}
0\lar A \xla{d} \fg\Tens A \xla{d} \textstyle\bigwedge^2\fg\Tens A \xla{d} \cdots
\textstyle\bigwedge^{p-1}\fg\Tens A \xla{d} \textstyle\bigwedge^p\fg\Tens A \xla{d}\cdots
\end{equation*}
The differential $d$ is given by
\begin{multline}
\label{stand_chain}
d(X_1\wedge\cdots\wedge X_p\otimes a)=
\sum_{i=1}^p (-1)^{i-1} X_1\wedge\cdots\wedge\hat X_i\wedge
\cdots\wedge X_p\otimes X_i\cdot a \\
+\sum_{1\le i<j\le p} (-1)^{i+j} [X_i,X_j]\wedge X_1\wedge\cdots\wedge\hat X_i
\wedge\cdots\wedge\hat X_j\wedge\cdots\wedge X_p\otimes a.
\end{multline}
(Here, as usual, the notation $\hat X_i$ indicates that $X_i$ is omitted.)
We also consider $\CC$ as a trivial $\fh$-module, and denote
the differential in the standard complex $C_\cdot(\fh,\CC)$ by $d'$.

In order to prove that $H_k(\fg,A)\ne 0$, it suffices to find a $k$-cycle
$z\in C_k(\fg,A)$ which is not a boundary. Note that
$\bigwedge^k\fg=\bigwedge^k\fh\oplus E$, where
\begin{equation*}
E=\bigoplus_{i=1}^k \textstyle\bigwedge^i\fr\Tens\textstyle\bigwedge^{k-i}\fh.
\end{equation*}
Fix $\eta\in\bigwedge^k\fh,\; \eta\ne 0$, and set $z=\eta\otimes 1\in C_k(\fg,A)$.
Since $\fh$ acts on $A$ by derivations, we have $X\cdot 1=0$ for each $X\in\fh$.
Now it follows from \eqref{stand_chain} that $d(\eta\otimes 1)=(d'\eta)\otimes 1$,
i.e., only the second group of summands in \eqref{stand_chain}
survives. On the other hand, since $\fh$ is semisimple, we have $H_k(\fh,\CC)\ne 0$,
i.e., the differential $d'\colon\bigwedge^k\fh\to\bigwedge^{k-1}\fh$ is zero
\cite{Koszul}. Therefore $d(\eta\otimes 1)=(d'\eta)\otimes 1=0$.

In order to prove that $\eta\otimes 1$ is not a boundary,
note that $A$ has a canonical augmentation
$\eps_A\colon A\to\CC$ defined by $\eps_A=\eps i_1$,
where $\eps$ is the counit of $\wh{U}(\fg)$,
and $i_1\colon A\to\wh{U}(\fg)=A\Psmash\wh{U}(\fh),\; a\mapsto a\otimes 1$
is the canonical embedding. Clearly, the restriction of $\eps_A$ to
$U(\fr)$ is precisely the counit of $U(\fr)$.
Now take $\xi\in \bigl(\bigwedge^k\fg\bigr)^*$ such that $\xi(\eta)=1$ and $\xi|_E=0$.
We then have $(\xi\otimes\eps_A)(\eta\otimes 1)=1$.
Let us show that $\xi\otimes\eps_A$ vanishes on $\Im d$.
To this end, consider the decomposition
\begin{equation}
\label{decomp_k+1}
\textstyle\bigwedge^{k+1}\fg\Tens A=
\bigl(\fr\Tens\textstyle\bigwedge^k\fh\Tens A\bigr)\oplus (F\Tens A),
\end{equation}
where
\begin{equation*}
F=\bigoplus_{i=2}^{k+1} \textstyle\bigwedge^i\fr\Tens\textstyle\bigwedge^{k+1-i}\fh.
\end{equation*}
It follows from \eqref{stand_chain} that $d(F\Tens A)\subset E\Tens A$.
By the same formula,
for each $X\in\fr$ and each $a\in A$ we have
$d(X\otimes\eta\otimes a)=\eta\otimes Xa+w$ for some $w\in E\Tens A$.
Since $\xi|_E=0$, we see that $\xi\otimes\eps_A$ vanishes on $E\Tens A$.
On the other hand, we have $\eps_A(Xa)=\eps_{U(\fr)}(X)\eps_A(a)=0$,
and so $(\xi\otimes\eps_A)(\eta\otimes Xa)=0$.
Together with \eqref{decomp_k+1}, this implies that
$\xi\otimes\eps_A$ vanishes on $\Im d$, and so
$\eta\otimes 1\notin\Im d$. The rest is clear.
\end{proof}

Combining Lemma~\ref{lemma:Uproj}, Lemma~\ref{lemma:ne0},
and Proposition~\ref{prop:Tor}, we obtain the following.

\begin{theorem}
Let $\fg$ be a finite-dimensional Lie algebra
such that $\wh{U}(\fg)$ is stably flat over $U(\fg)$.
Then $\fg$ is solvable.
\end{theorem}

\vspace*{10mm}
\begin{flushleft}
\scshape\small
Department of Differential Equations and Functional Analysis\\
Faculty of Science\\
Peoples' Friendship University of Russia\\
Mikluho-Maklaya 6\\
117198 Moscow\\
RUSSIA

\medskip
{\itshape Address for correspondence:}\\

\medskip\upshape
Krupskoi 8--3--89\\
Moscow 119311\\
Russia

\medskip
{\itshape E-mail:} {\ttfamily pirkosha@sci.pfu.edu.ru, pirkosha@online.ru}
\end{flushleft}

\begin{thebibliography}{XX}
\bibitem{BFGP}
Bonneau, P., Flato, M., Gerstenhaber, M., Pinczon, G.
{\em The hidden group structure of quantum groups:
strong duality, rigidity and preferred deformations},
Comm. Math. Phys. \textbf{161} (1994), 125--156.
\bibitem{DKR}
Doplicher, S., Kastler, D., Robinson, D. W.
{\em Covariance algebras in field theory and statistical mechanics}.
Comm. Math. Phys. \textbf{3} (1966), 1--28.
\bibitem{Dos_hprop}
Dosiev, A. A.
{\em Homological dimensions of the algebra formed by entire
functions of elements of a nilpotent Lie algebra},
Funct. Anal. Appl. \textbf{37} (1), 61--64.
\bibitem{X1}
Helemskii, A. Ya. {\itshape The Homology of Banach and Topological Algebras},
Moscow University Press, 1986 (Russian); English transl.: Kluwer Academic
Publishers, Dordrecht, 1989.
\bibitem{X2}
Helemskii, A. Ya. {\itshape Banach and Polynormed Algebras: General Theory,
Representations, Homology}, Nauka, Moscow, 1989 (Russian); English transl.:
Oxford University Press, 1993.
\bibitem{KV}
Kiehl, R. and Verdier, J. L.
{\itshape Ein einfacher Beweis des Koh\"arenzsatzes von Grauert},
Math. Ann. \textbf{195} (1971), 24--50.
\bibitem{Koszul}
Koszul, J.-L.
{\em Homologie et cohomologie des alg\`ebres de Lie},
Bull. Soc. Math. France \textbf{78} (1950), 65--127.
\bibitem{MR}
McConnell, J. C.; Robson, J. C.
{\em Noncommutative Noetherian rings}.
John Wiley \& Sons, Ltd., Chichester, 1987.
\bibitem{MRL}
Mitiagin, B.; Rolewicz, S.; \.Zelazko, W.
{\em Entire functions in $B_0$-algebras.}
Studia Math. \textbf{21} (1961/1962), 291--306.
\bibitem{NR}
Neeman, A. and Ranicki, A.
{\em Noncommutative localization and chain complexes. I.
Algebraic $K$- and $L$-theory},
Preprint arXiv.org:math.RA/0109118.
\bibitem{Pir_stb}
Pirkovskii, A. Yu.
{\em Stably flat completions of universal enveloping algebras},
Preprint arXiv.org:math.FA/0311492.
\bibitem{Schwtzr}
Schweitzer, L. B.
{\em Dense $m$-convex Fr\'echet subalgebras of operator algebra
crossed products by Lie groups}.
Internat. J. Math. \textbf{4} (4) (1993), 601--673.
\bibitem{Sweedler}
Sweedler, M. E.
{\em Hopf algebras}.
Benjamin, New York, 1969.
\bibitem{T1}
Taylor, J. L. {\itshape Homology and cohomology for topological algebras},
Adv. Math. \textbf{9} (1972), 137--182.
\bibitem{T2}
Taylor, J. L. {\itshape A general framework for a multi-operator functional calculus},
Adv. Math. \textbf{9} (1972), 183--252.
\bibitem{T3}
Taylor, J. L.
{\itshape Functions of several noncommuting variables},
Bull. Amer. Math. Soc. \textbf{79} (1973), 1--34.
\end{thebibliography}
\end{document}